# MAXIMAL SOLVABLE LEIBNIZ ALGEBRAS WHOSE NILRADICAL IS A QUASI-FILIFORM ALGEBRA

ABDURASULOV K.K.[1,2] AND ADASHEV J.Q.[1,2,3]

[1]*Institute of Mathematics, Uzbekistan Academy of Sciences, 4 University street, 100174.*
abdurasulov0505@mail.ru, adashevjq@mail.ru

[2] *National University of Uzbekistan, 4 University street, 100174, Tashkent, Uzbekistan.*

[3]*Chirchiq State Pedagogical Institute of Tashkent region, 104 Amir Temur street, 111700, Tashkent, Uzbekistan.*

ABSTRACT. The present article is a part of the study of solvable Leibniz algebras with a given nilradical. In this paper solvable Leibniz algebras, whose nilradicals is naturally graded quasi-filiform algebra and the complemented space to the nilradical has maximal dimension, are described up to isomorphism.

## 1. INTRODUCTION

During the last decades the theory of Leibniz algebras has been actively investigated and many results of the Lie Theory have been transferred to Leibniz algebras.

Levi's decomposition asserts that every finite-dimensional Lie algebra is a semidirect sum of a semisimple Lie subalgebra and solvable radical, while semisimple Lie algebras over the field of complex numbers have been classified by E. Cartan [16] and over the field of real numbers by F. Gantmacher [19]. Thus, the problem of a description of finite-dimensional Lie algebras is reduced to the study of solvable Lie algebras. At the same period, an essential progress was made by Mal'cev reducing the problem of classification of solvable Lie algebras to nilpotent Lie algebras. Since then classification results have been all related to the nilpotent part.

Using a result of [27], an approach to the study of solvable Lie algebras in an arbitrary finite-dimensions through the use of the nilradical was developed in [4, 28, 29, 35, 36], etc. In particular, García, L. [20] studied solvable Lie algebras with quasi-filiform nilradicals. In fact, there are solvable Lie algebras constructed using the method explained in [27].

Leibniz algebras, a "noncommutative version" of Lie algebras, were first introduced in 1965 by A. Bloh [8] under the name "$D$-algebras". They appeared again in 1993 after Loday's work [26], where he reintroduced them, coining the term "Leibniz algebra". Recently, it has been a trend to show how various results from Lie algebras extend to Leibniz algebras. In particular, it has been interest in extending classifications of certain classes of Lie algebras to classifications of corresponding Leibniz algebras [6, 7, 9–12, 15, 25]. For finite-dimensional Leibniz algebras over a field of characteristic zero, there is an analogue of Levi's decomposition: namely, any Leibniz algebra is decomposed into a semidirect sum of a semisimple Lie algebra and its solvable radical [7]. Therefore, similar to Lie case, the main problem of the study of Leibniz algebras is reduced to solvable ones.

The analogue of Mubarakzjanov's result has been applied to the Leibniz algebras in [17], showing the importance of consideration of the nilradical in the case of Leibniz algebras as well. Papers [1–3, 9, 13, 14, 17, 18, 23, 24, 31–34] are also devoted to the study of solvable Leibniz algebras by considering the nilradical.

The aim of this article is to describe solvable Leibniz algebras with naturally graded quasi-filiform Leibniz nilradicals and with maximal dimension of complemented space of its nilradical. Namely, naturally graded quasi-filiform Leibniz algebras in any finite dimension over $\mathbb{C}$ were studied by Camacho, Gómez, González, and Omirov [11]. They found five such algebras of the first type, where two of them depend on a parameter and eight algebras of the second type with one of them depending on a parameter. The naturally graded quasi-filiform Lie algebras were classified in [21]. Here exist six families, two of them are decomposable, i.e., split into a direct sum of ideals and as well as there exist some special cases that appear only in low dimensions.

It is know that in works devoted to the classification of solvable Leibnits algebras generated by their nilradicals, algebras with certain nilradical have been studied. In this work, algebras that their nilradicals are isomorphic to quasi-filiform algebras are studied. It should be noted that the previous result was used directly to obtain the solvable algebra, so the computational processes were much simpler than in the previous work.

Throughout the paper vector spaces and algebras are finite-dimensional over the field of the complex numbers. Moreover, in the table of multiplication of any algebra the omitted products are assumed to be zero and, if it is not noted, we consider non-nilpotent solvable algebras.

## 2. Preliminaries

In this section we recall some basic notions and concepts used throughout the paper.

**Definition 2.1.** *A vector space with a bilinear bracket $(L, [\cdot, \cdot])$ is called a Leibniz algebra if for any $x, y, z \in L$ the so-called Leibniz identity*
$$[x, [y, z]] = [[x, y], z] - [[x, z], y],$$
*holds.*

Further we use the notation
$$\mathcal{L}I(x, y, z) = [x, [y, z]] - [[x, y], z] + [[x, z], y].$$
It is obvious that Leibniz algebras are determined by the identity $\mathcal{L}I(x, y, z) = 0$.

The sets $Ann_r(L) = \{x \in L : [y, x] = 0, \ \forall y \in L\}$ and $Ann_l(L) = \{x \in L : [x, y] = 0, \ \forall y \in L\}$ are called the right and left annihilators of $L$, respectively. It is observed that $Ann_r(L)$ is a two-sided ideal of $L$, and for any $x, y \in L$ the elements $[x, x]$ and $[x, y] + [y, x]$ belong to $Ann_r(L)$.

For a given Leibniz algebra $(L, [\cdot, \cdot])$, the sequences of two-sided ideals are defined recursively as follows:
$$L^1 = L, \ L^{k+1} = [L^k, L], \ k \geq 1, \qquad L^{[1]} = L, \ L^{[s+1]} = [L^{[s]}, L^{[s]}], \ s \geq 1.$$

These are said to be the lower central and the derived series of $L$, respectively.

**Definition 2.2.** *A Leibniz algebra $L$ is said to be nilpotent (respectively, solvable), if there exists $n \in \mathbb{N}$ ($m \in \mathbb{N}$) such that $L^n = \{0\}$ (respectively, $L^{[m]} = \{0\}$).*

It is easy to see that the sum of two nilpotent ideals is nilpotent. Therefore, the maximal nilpotent ideal always exists. The maximal nilpotent ideal of a Leibniz algebra is said to be the nilradical of the algebra.

**Definition 2.3.** *A linear map $d \colon L \to L$ of a Leibniz algebra $(L, [\cdot, \cdot])$ is said to be a derivation if for any $x, y \in L$, the following condition holds:*
$$d([x, y]) = [d(x), y] + [x, d(y)]. \tag{2.1}$$

The set of all derivations of $L$ is denoted by $\text{Der}(L)$, which is a Lie algebra with respect to the commutator.

For a given element $x$ of a Leibniz algebra $L$, the right multiplication operator $\mathcal{R}_x \colon L \to L$, defined by $\mathcal{R}_x(y) = [y, x], y \in L$ is a derivation. In fact, Leibniz algebras are characterized by this property regarding right multiplication operators. As in the Lie case, these kinds of derivations are said to be inner derivations.

**Definition 2.4.** *Let $d_1, d_2, \ldots, d_n$ be derivations of a Leibniz algebra $L$. The derivations $d_1, d_2, \ldots, d_n$ are said to be linearly nil-independent, if for $\alpha_1, \alpha_2, \ldots, \alpha_n \in \mathbb{C}$ and a natural number $k$,*
$$(\alpha_1 d_1 + \alpha_2 d_2 + \cdots + \alpha_n d_n)^k = 0 \ \text{implies} \ \alpha_1 = \alpha_2 = \cdots = \alpha_n = 0.$$

Note that in the above definition the power is understood with respect to composition.

Let $L$ be a solvable Leibniz algebra. Then it can be written in the form $L = N \oplus Q$, where $N$ is the nilradical and $Q$ is the complementary subspace.

**Theorem 2.5** ( [17]). *Let $L$ be a solvable Leibniz algebra and $N$ be its nilradical. Then the dimension of $Q$ is not greater than the maximal number of nil-independent derivations of $N$.*

In other words, similar to the Lie algebras we have that $dim L \leq 2 dim N$. Note that the quotient Leibniz algebra $L/N^2$ is solvable, with the abelian nilradical $N/N^2$. Thus we have that $L/N^2$ can be written as a vector space $L/N^2 = \overline{Q} \oplus N/N^2$, with $dim\overline{Q} \leq dim(N/N^2)$. Moreover all generator elements of $N$ belong to $N \setminus N^2$, we can conclude that linear nil-independent derivations of $N$ induced nil-independent derivations of $N/N^2$. Thus we have $dimQ \leq dim\overline{Q}$. Therefore, for the solvable algebra $L = N \oplus Q$ we obtain that $dimQ \leq dim(N/N^2)$.

Below we define the notion of a quasi-filiform Leibniz algebra.

**Definition 2.6.** *A Leibniz algebra $L$ is called quasi-filiform if $L^{n-2} \neq \{0\}$ and $L^{n-1} = \{0\}$, where $n = \dim L$.*

Given an $n$-dimensional nilpotent Leibniz algebra $L$ such that $L^{s-1} \neq \{0\}$ and $L^s = \{0\}$, put $L_i = L^i/L^{i+1}$, $1 \leq i \leq s-1$, and $gr(L) = L_1 \oplus L_2 \oplus \cdots \oplus L_{s-1}$. Using $[L^i, L^j] \subseteq L^{i+j}$, it is easy to establish that $[L_i, L_j] \subseteq L_{i+j}$. So, we obtain the graded algebra $gr(L)$. If $gr(L)$ and $L$ are isomorphic, then we say that $L$ is *naturally graded*.

Let $x$ be a nilpotent element of the set $L \setminus L^2$. For the nilpotent operator of right multiplication $\mathcal{R}_x$ we define a decreasing sequence $C(x) = (n_1, n_2, \ldots, n_k)$, where $n = n_1 + n_2 + \cdots + n_k$, which consists of the dimensions of Jordan blocks of the operator $\mathcal{R}_x$. In the set of such sequences we consider the lexicographic order, that is, $C(x) = (n_1, n_2, \ldots, n_k) \leq C(y) = (m_1, m_2, \ldots, m_t) \Leftrightarrow$ there exists $i \in \mathbb{N}$ such that $n_j = m_j$ for any $j < i$ and $n_i < m_i$.

**Definition 2.7.** *The sequence $C(L) = \max\limits_{x \in L \setminus L^2} C(x)$ is called a characteristic sequence of the algebra $L$.*

Let $L$ be an $n$-dimensional naturally graded quasi-filiform non-Lie Leibniz algebra which has the characteristic sequence $(n-2, 1, 1)$ or $(n-2, 2)$. The first case (case 2-filiform) has been studied in [12] and the second in [11]. Thanks to [3], we already have the classification of solvable Leibniz algebras whose nilradical is a 2-filiform Leibniz algebra.

**Definition 2.8.** *A quasi-filiform Leibniz algebra $L$ is called algebra of the type I (respectively, type II), if there exists a basic element $e_1 \in L \setminus L^2$ such that the operator $\mathcal{R}_{e_1}$ has the form $\begin{pmatrix} J_{n-2} & 0 \\ 0 & J_2 \end{pmatrix}$ (respectively, $\begin{pmatrix} J_2 & 0 \\ 0 & J_{n-2} \end{pmatrix}$).*

In the following theorem we give the classification of naturally graded quasi-filiform Leibniz algebras given in [11].

**Theorem 2.9.** *An arbitrary n-dimensional naturally graded quasi-filiform Leibniz algebra of type I is isomorphic to one of the following pairwise non-isomorphic algebras of the families:*

$$\mathcal{L}_n^{1,\beta} : \begin{cases} [e_i, e_1] = e_{i+1}, \ 1 \leq i \leq n-3 \\ [e_{n-1}, e_1] = e_n, \\ [e_1, e_{n-1}] = \beta e_n, \ \beta \in \mathbb{C} \end{cases} \qquad \mathcal{L}_n^{2,\beta} : \begin{cases} [e_i, e_1] = e_{i+1}, \ 1 \leq i \leq n-3 \\ [e_{n-1}, e_1] = e_n, \\ [e_1, e_{n-1}] = \beta e_n, \ \beta \in \{0, 1\} \\ [e_{n-1}, e_{n-1}] = e_n \end{cases}$$

$$\mathcal{L}_n^{3,\beta} : \begin{cases} [e_i, e_1] = e_{i+1}, \ 1 \leq i \leq n-3 \\ [e_{n-1}, e_1] = e_n + e_2, \\ [e_1, e_{n-1}] = \beta e_n, \ \beta \in \{-1, 0, 1\} \end{cases} \qquad \mathcal{L}_n^{4,\gamma} : \begin{cases} [e_i, e_1] = e_{i+1}, \ 1 \leq i \leq n-3 \\ [e_{n-1}, e_1] = e_n + e_2, \\ [e_{n-1}, e_{n-1}] = \gamma e_n, \ \gamma \neq 0 \end{cases}$$

$$\mathcal{L}_n^{5,\beta,\gamma} : \begin{cases} [e_i, e_1] = e_{i+1}, \ 1 \leq i \leq n-3 \\ [e_{n-1}, e_1] = e_n + e_2, \\ [e_1, e_{n-1}] = \beta e_n, \ (\beta, \gamma) = (1,1) \text{ or } (2,4) \\ [e_{n-1}, e_{n-1}] = \gamma e_n, \end{cases}$$

*where $\{e_1, e_2, \ldots, e_n\}$ is a basis of the algebra.*

**Theorem 2.10.** *An arbitrary n-dimensional naturally graded quasi-filiform Leibniz algebra of second II type is isomorphic to one of the following pairwise non-isomorphic algebras of the families:*

$n$ *even*

$$\mathcal{L}_n^1 : \begin{cases} [e_1, e_1] = e_2, \\ [e_i, e_1] = e_{i+1},\ 3 \leq i \leq n-1, \\ [e_1, e_i] = -e_{i+1},\ 3 \leq i \leq n-1, \end{cases} \qquad \mathcal{L}_n^2 : \begin{cases} [e_1, e_1] = e_2, \\ [e_i, e_1] = e_{i+1},\ 3 \leq i \leq n-1, \\ [e_1, e_3] = e_2 - e_4, \\ [e_1, e_i] = -e_{i+1},\ 4 \leq i \leq n-1, \end{cases}$$

$$\mathcal{L}_n^3 : \begin{cases} [e_1, e_1] = e_2, \\ [e_i, e_1] = e_{i+1},\ 3 \leq i \leq n-1, \\ [e_1, e_i] = -e_{i+1},\ 3 \leq i \leq n-1, \\ [e_3, e_3] = e_2, \end{cases} \qquad \mathcal{L}_n^4 : \begin{cases} [e_1, e_1] = e_1, \\ [e_i, e_1] = e_{i+1},\ 3 \leq i \leq n-1, \\ [e_1, e_3] = 2e_2 - e_4, \\ [e_1, e_i] = -e_{i+1},\ 4 \leq i \leq n-1, \\ [e_3, e_3] = e_2, \end{cases}$$

$n$ odd, $\mathcal{L}_n^1, \mathcal{L}_n^2, \mathcal{L}_n^3, \mathcal{L}_n^4$,

$$\mathcal{L}_n^5 : \begin{cases} [e_1, e_1] = e_2, \\ [e_i, e_1] = e_{i+1},\ 3 \leq i \leq n-1, \\ [e_1, e_i] = -e_{i+1},\ 3 \leq i \leq n-1, \\ [e_i, e_{n+2-i}] = (-1)^i e_n,\ 3 \leq i \leq n-1, \end{cases} \qquad \mathcal{L}_n^{6,\beta} : \begin{cases} [e_1, e_1] = e_2, \\ [e_i, e_1] = e_{i+1},\ 3 \leq i \leq n-1, \\ [e_1, e_3] = \beta e_2 - e_4,\ \beta \in \{1, 2\}, \\ [e_1, e_i] = -e_{i+1},\ 4 \leq i \leq n-1, \\ [e_i, e_{n+2-i}] = (-1)^i e_n,\ 3 \leq i \leq n-1, \end{cases}$$

$$\mathcal{L}_n^{7,\gamma} : \begin{cases} [e_1, e_1] = e_2, \\ [e_i, e_1] = e_{i+1},\ 3 \leq i \leq n-1, \\ [e_1, e_i] = -e_{i+1},\ 3 \leq i \leq n-1, \\ [e_3, e_3] = \gamma e_2,\ \gamma \neq 0, \\ [e_i, e_{n+2-i}] = (-1)^i e_n,\ 3 \leq i \leq n-1, \end{cases} \qquad \mathcal{L}_n^{8,\beta,\gamma} : \begin{cases} [e_1, e_1] = e_2, \\ [e_i, e_1] = e_{i+1},\ 3 \leq i \leq n-1, \\ [e_1, e_3] = \beta e_2 - e_4, \\ [e_1, e_i] = -e_{i+1},\ 4 \leq i \leq n-1, \\ [e_3, e_3] = \gamma e_2,\ (\beta, \gamma) = (-2, 1), (2, 1)\ or\ (4, 2), \\ [e_i, e_{n+2-i}] = (-1)^i e_n,\ 3 \leq i \leq n-1, \end{cases}$$

where $\{e_1, e_2, \ldots, e_n\}$ is a basis of the algebra.

The study of naturally graded quasi-filiform Leibniz algebra of corresponding type in Theorems 2.9 and 2.10 can be simplified, as follows (see [10]):

**Proposition 2.11.** *Let $L$ be a naturally graded quasi-filiform Leibniz algebra, then it is isomorphic to one algebra of the non isomorphic families*

$$\mathcal{L}(\alpha, \beta, \gamma) : \begin{cases} [e_i, e_1] = e_{i+1}, & 1 \leq i \leq n-3, \\ [e_{n-1}, e_1] = e_n + \alpha e_2, & [e_1, e_{n-1}] = \beta e_n,\ [e_{n-1}, e_{n-1}] = \gamma e_n, \end{cases}$$

$$\mathcal{G}(\alpha, \beta, \gamma) : \begin{cases} [e_1, e_1] = e_2, & [e_i, e_1] = e_{i+1}, & 3 \leq i \leq n-1, \\ [e_1, e_3] = -e_4 + \beta e_2, & [e_1, e_i] = -e_{i+1}, & 4 \leq i \leq n-1, \\ [e_3, e_3] = \gamma e_2, & [e_i, e_{n+2-i}] = (-1)^i \alpha e_n, & 3 \leq i \leq n-1, \end{cases}$$

*where $\{e_1, e_2, \ldots, e_n\}$ is a basis of the algebra and in the algebra $\mathcal{G}(\alpha, \beta, \gamma)$ if $n$ is odd, then $\alpha \in \{0, 1\}$, if $n$ is even, then $\alpha = 0$.*

**Remark 2.12.** *The algebras given in Theorem 2.9 and 2.10 which stated in Proposition 2.11 are of the form:*

$\mathcal{L}(0, \beta, 0) := \mathcal{L}_n^{1,\beta}; \quad \mathcal{L}(0, \beta, 1) := \mathcal{L}_n^{2,\beta}; \quad \mathcal{L}(1, \beta, 0) := \mathcal{L}_n^{3,\beta}; \quad \mathcal{L}(1, 0, \gamma) := \mathcal{L}_n^{4,\gamma}; \quad \mathcal{L}(1, \beta, \gamma) := \mathcal{L}_n^{5,\beta,\gamma};$

$\mathcal{G}(0, 0, 0) := \mathcal{L}_n^1; \quad \mathcal{G}(0, 1, 0) := \mathcal{L}_n^2; \quad \mathcal{G}(0, 0, 1) := \mathcal{L}_n^3; \quad \mathcal{G}(0, 2, 1) := \mathcal{L}_n^4; \quad \mathcal{G}(1, 0, 0) := \mathcal{L}_n^5;$

$\mathcal{G}(1, \beta, 0) := \mathcal{L}_n^{6,\beta}; \quad \mathcal{G}(1, 0, \gamma) := \mathcal{L}_n^{7,\gamma}; \quad \mathcal{G}(1, \beta, \gamma) := \mathcal{L}_n^{8,\beta,\gamma}.$

## 3. Main results. Solvable Leibniz algebras whose nilradical is quasi-filiform non-Lie Leibniz algebras.

This section is devoted to the classification of solvable Leibniz algebras whose nilradical is naturally graded quasi-filiform Leibniz algebras. Due to Proposition 2.11 we only need to consider solvable Leibniz algebras with nilradicals $\mathcal{L}(\alpha, \beta, \gamma)$ and $\mathcal{G}(\alpha, \beta, \gamma)$.

### 3.1. Derivations of algebras $\mathcal{L}(\alpha, \beta, \gamma)$ and $\mathcal{G}(\alpha, \beta, \gamma)$.

In order to start the description we need to know the derivations of naturally graded quasi-filiform Leibniz algebras.

**Proposition 3.1.** *An arbitrary $d \in Der(\mathcal{L}(\alpha, \beta, \gamma))$ has the following form:*

$$\begin{cases} d(e_1) = \sum_{t=1}^{n} a_t e_t, \\ d(e_2) = (2a_1 + a_{n-1}\alpha)e_2 + \sum_{t=3}^{n-2} a_{t-1}e_t + (a_{n-1} + a_{n-1}\beta)e_n, \\ d(e_i) = (ia_1 + a_{n-1}\alpha)e_i + \sum_{t=i+1}^{n-2} a_{t-i+1}e_t, \quad 3 \leq i \leq n-2, \\ d(e_{n-1}) = \sum_{t=2}^{n} b_t e_t, \\ d(e_n) = (b_{n-3} - a_{n-3}\alpha)e_{n-2} + (b_{n-1} + a_1 + a_{n-1}\gamma - a_{n-1}\alpha(1+\beta))e_n, \end{cases}$$

*where*

$$\begin{aligned} &b_i = a_i\alpha, \ 2 \leq i \leq n-4, \ \beta(b_{n-3} - a_{n-3}\alpha) = \gamma(b_{n-3} - a_{n-3}\alpha) = 0, \ b_{n-1}\alpha = a_1\alpha + a_{n-1}\alpha^2, \\ &\gamma b_{n-1} = \gamma(a_1 + a_{n-1}\gamma - a_{n-1}\alpha(1+\beta)), \quad \gamma a_{n-1} = \beta a_{n-1}(\gamma - \alpha(1+\beta)). \end{aligned} \quad (3.1)$$

*Proof.* It is easy to see that $\{e_1, e_{n-1}\}$ are the generator basis elements of the algebra $\mathcal{L}(\alpha, \beta, \gamma)$.
We put

$$d(e_1) = \sum_{t=1}^{n} a_t e_t, \qquad d(e_{n-1}) = \sum_{t=1}^{n} b_t e_t.$$

From the derivation property (2.1) we have

$$d(e_2) = d([e_1, e_1]) = [d(e_1), e_1] + [e_1, d(e_1)] = (2a_1 + a_{n-1}\alpha)e_2 + \sum_{t=3}^{n-2} a_{t-1}e_t + (a_{n-1} + a_{n-1}\beta)e_n.$$

By the induction and the property of derivation (2.1) we derive

$$d(e_i) = (ia_1 + a_{n-1}\alpha)e_i + \sum_{t=i+1}^{n-2} a_{t-i+1}e_t, \quad 3 \leq i \leq n-2.$$

From the derivation property (2.1) we have

$$d(e_n) = d([e_{n-1}, e_1]) - \alpha d(e_2) =$$

$$= (b_1 + b_{n-1}\alpha - a_1\alpha - a_{n-1}\alpha^2)e_2 + \sum_{t=3}^{n-2}(b_{t-1} - a_{t-1}\alpha)e_t + (b_{n-1} + a_1 + a_{n-1}\gamma - a_{n-1}\alpha(1+\beta))e_n.$$

Considering

$$0 = d([e_2, e_{n-1}]) = d([e_n, e_1]).$$

Consequently,

$$b_1 = 0, \ b_{n-1}\alpha = a_1\alpha + a_{n-1}\alpha^2, \ b_i = a_i\alpha, \ 2 \leq i \leq n-4.$$

Using property of the derivation for the products $[e_1, e_{n-1}] = \beta e_n$, $[e_{n-1}, e_{n-1}] = \gamma e_n$, we have

$$\beta(b_{n-3} - a_{n-3}\alpha) = \gamma(b_{n-3} - a_{n-3}\alpha) = 0,$$
$$\gamma b_{n-1} = \gamma(a_1 + a_{n-1}\gamma - a_{n-1}\alpha(1+\beta)), \quad \gamma a_{n-1} = \beta a_{n-1}(\gamma - \alpha(1+\beta)).$$

$\square$

**Proposition 3.2.** *Any derivation of the algebras $\mathcal{G}(\alpha, \beta, \gamma)$ has the following form:*

$$\begin{cases} d(e_1) = \sum_{t=1}^{n} a_t e_t, \\ d(e_3) = \sum_{t=2}^{n} b_t e_t, \\ d(e_2) = (2a_1 + a_3\beta)e_2, \\ d(e_4) = \gamma a_3 e_2 + (a_1 + b_3)e_4 + \sum_{t=5}^{n-1} b_{t-1}e_t + (b_{n-1} - a_{n-1}\alpha)e_n, \\ d(e_i) = ((i-3)a_1 + b_3)e_i + \sum_{t=i+1}^{n-1} b_{t-i+3}e_t + (b_{n-i+3} - (-1)^i a_{n-i+3}\alpha)e_n, \ 5 \leq i \leq n-1, \\ d(e_n) = ((n-3)a_1 + b_3 - (-1)^n a_3\alpha)e_n, \end{cases}$$

*where*

$$2a_3\gamma+b_3\beta = a_1\beta+a_3\beta^2, \quad (1-(-1)^n)a_{n-1}\alpha = 0, \quad 2b_3\gamma = \gamma(2a_1+a_3\beta), \quad b_3\alpha = a_1\alpha-(-1)^n a_3\alpha^2. \quad (3.2)$$

*Proof.* From Proposition 2.11 we conclude that $e_1$ and $e_3$ are the generator basis elements of the algebra.

We put

$$d(e_1) = \sum_{t=1}^{n} a_t e_t,, \qquad d(e_3) = \sum_{t=1}^{n} b_t e_t.$$

From the derivation property (2.1) we have

$$d(e_2) = d([e_1, e_1]) = (2a_1 + a_3\beta)e_2,$$

$$d(e_4) = d([e_3, e_1]) = [d(e_3), e_1] + [e_3, d(e_1)] =$$

$$= [\sum_{t=1}^{n} b_t e_t, e_1] + [e_3, \sum_{t=1}^{n} a_t e_t] = b_1 e_2 + \sum_{t=4}^{n} b_{t-1} e_t + a_1 e_4 + \gamma a_3 e_2 - a_{n-1}\alpha e_n$$

$$= (a_1 + b_3)e_4 + \sum_{t=5}^{n-1} b_{t-1} e_t + (b_1 + \gamma a_3)e_2 + (b_{n-1} - a_{n-1}\alpha)e_n.$$

Applying the induction and the derivation property, we derive

$$d(e_i) = ((i-3)a_1 + b_3)e_i + \sum_{t=i+1}^{n-1} b_{t-i+3} e_t + (b_{n-i+3} - (-1)^i a_{n-i+3}\alpha)e_n, \quad 4 \leq i \leq n-1,$$

$$d(e_n) = ((n-3)a_1 + b_3 - (-1)^n a_3\alpha)e_n.$$

From $0 = d([e_i, e_3]) = [d(e_i), e_3] + [e_i, d(e_3)]$, $3 \leq i \leq n-2$, we conclude

$$b_1 = 0, \quad ((-1)^i - 1)\alpha b_{n-i+2} = 0, \quad 3 \leq i \leq n-2.$$

Using the derivation for the products $[e_1, e_3] = -e_4 + \beta e_2$, $[e_3, e_3] = \gamma e_2$, $[e_3, e_{n-1}] = -\alpha e_n$, we have

$$2a_3\gamma + b_3\beta = a_1\beta + a_3\beta^2, \quad (1-(-1)^n)a_{n-1}\alpha = 0, \quad 2b_3\gamma = \gamma(2a_1 + a_3\beta), \quad b_3\alpha = a_1\alpha - (-1)^n a_3\alpha^2.$$

□

The following theorem describes the maximal dimensions of the complemented spaces to $\mathcal{L}(\alpha, \beta, \gamma)$ and $\mathcal{G}(\alpha, \beta, \gamma)$.

**Theorem 3.3.** *Let $R$ be a solvable Leibniz algebra whose nilradical is naturally graded quasi-filiform non-Lie Leibniz algebras. Then the maximal dimension of complemented space to the nilradical are not greater than two.*

*Proof.* Due to Propositions 3.1-3.2 the nilpotency of derivation of naturally graded quasi-filiform non-Lie Leibniz algebras depends on the following parametrs:

- For $\mathcal{L}(\alpha, \beta, \gamma)$, the nilpotency of derivation depends on $a_1$ and $b_{n-1}$, i.e. the derivation is nilpotent if and only if $a_1 = b_{n-1} = 0$.
- For $\mathcal{G}(\alpha, \beta, \gamma)$, the nilpotency of derivation depends on $a_1$ and $b_3$, i.e. the derivation is nilpotent if and only if $a_1 = b_3 = 0$.

Applying the Theorem 2.5, the stated inequalities follow. □

**Remark 3.4.** *From the equations (3.1)-(3.2) and using the Theorem 3.3 it follows that for the possible values of the parameters $\alpha, \beta$ and $\gamma$, we derive the following table:*

TABLE 1. The dimensions of the complemented spaces to $\mathcal{L}(\alpha,\beta,\gamma)$ and $\mathcal{G}(\alpha,\beta,\gamma)$.

| Algebra | restrictions | dimensional of complementary space |
|---|---|---|
| $\mathcal{L}(0,\beta,0)$ | $b_i = 0,\ 2 \leq i \leq n-4,\ \beta b_{n-3} = 0,$ | $dimQ \leq 2,$ |
| $\mathcal{L}(0,0,1)$ | $a_{n-1} = b_i = 0,\ 2 \leq i \leq n-3,\ b_{n-1} = a_1,$ | $dimQ = 1,$ |
| $\mathcal{L}(0,1,1)$ | $b_i = 0,\ 2 \leq i \leq n-3,\ b_{n-1} = a_1 + a_{n-1},$ | $dimQ \leq 2,$ |
| $\mathcal{L}(1,-1,0)$ | $b_i = a_i,\ 2 \leq i \leq n-3,\ b_{n-1} = a_1 + a_{n-1},$ | $dimQ \leq 2,$ |
| $\mathcal{L}(1,0,0)$ | $b_i = a_i,\ 2 \leq i \leq n-4,\ b_{n-1} = a_1 + a_{n-1},$ | $dimQ \leq 2,$ |
| $\mathcal{L}(1,1,0)$ | $a_{n-1} = 0,\ b_i = a_i,\ 2 \leq i \leq n-3,\ b_{n-1} = a_1,$ | $dimQ = 1,$ |
| $\mathcal{L}(1,0,\gamma),\ \gamma \neq 0$ | $a_{n-1} = 0,\ b_i = a_i,\ 2 \leq i \leq n-3,\ b_{n-1} = a_1,$ | $dimQ = 1,$ |
| $\mathcal{L}(1,1,1)$ | $a_{n-1} = 0,\ b_i = a_i,\ 2 \leq i \leq n-3,\ b_{n-1} = a_1,$ | $dimQ = 1,$ |
| $\mathcal{L}(1,2,4)$ | $a_{n-1} = 0,\ b_i = a_i,\ 2 \leq i \leq n-3,\ b_{n-1} = a_1,$ | $dimQ = 1,$ |
| $\mathcal{G}(0,0,0)$ |  | $dimQ \leq 2,$ |
| $\mathcal{G}(0,1,0)$ | $b_3 = a_1 + a_3,$ | $dimQ \leq 2,$ |
| $\mathcal{G}(0,0,1)$ | $b_3 = a_1,\ a_3 = 0,$ | $dimQ = 1,$ |
| $\mathcal{G}(0,2,1)$ | $b_3 = a_1 + a_3,$ | $dimQ \leq 2,$ |
| $\mathcal{G}(1,0,0)$ | $b_3 = a_1 + a_3,\ a_{n-1} = 0,$ | $dimQ \leq 2,$ |
| $\mathcal{G}(1,1,0)$ | $b_3 = a_1 + a_3,\ a_{n-1} = 0,$ | $dimQ \leq 2,$ |
| $\mathcal{G}(1,2,0)$ | $b_3 = a_1,\ a_3 = 0,\ a_{n-1} = 0,$ | $dimQ = 1,$ |
| $\mathcal{G}(1,0,\gamma),\ \gamma \neq 0$ | $b_3 = a_1,\ a_3 = 0,\ a_{n-1} = 0,$ | $dimQ = 1,$ |
| $\mathcal{G}(1,-2,1)$ | $b_3 = a_1,\ a_3 = 0,\ a_{n-1} = 0,$ | $dimQ = 1,$ |
| $\mathcal{G}(1,2,1)$ | $b_3 = a_1 + a_3,\ a_{n-1} = 0,$ | $dimQ \leq 2,$ |
| $\mathcal{G}(1,4,2)$ | $b_3 = a_1, a_3 = 0,\ a_{n-1} = 0,$ | $dimQ = 1,$ |

3.2. **Solvable Leibniz algebras with codimensional nilradical equal to the number of generators nilradicals.**

We give a description of solvable Leibniz algebras such that the dimension of the complementary subspace is equal to the number of generators of nilradical. In other words, we describe solvable Leibniz algebras $R = N \oplus Q$ with $\dim Q = \dim N/N^2$.

Let the multiplication table of the nilradical $N$ be expressed through the products:

$$[e_i, e_j] = \sum_{t=k+1}^{n} \gamma_{i,j}^t e_t, \quad 1 \leq i, j \leq n.$$

Let $\{e_{i_1}, e_{i_2}, \ldots, e_{i_{n_i}}\}$ be a basis of the space $N_i := span(N^i \setminus N^{i+1})$, $1 \leq i \leq s-1$ and $\dim N_i = n_i$, where $n_1 = k$ and $n_1 + n_2 + \cdots + n_{s-1} = n$. Now we give a description of the solvable Leibniz algebras with a codimensional nilradical equal to the number of generator basis elements of nilradicals.

**Theorem 3.5.** *[1] Let $R = N \oplus Q$ be a solvable Leibniz algebra such that $dimQ = dimN/N^2 = k$. Then $R$ admits a basis $\{e_1, e_2, \ldots, e_n, x_1, \ldots, x_k\}$ such that the table of multiplication in $R$ has the following form:*

$$\begin{cases} [e_i, e_j] = \sum_{t=k+1}^{n} \gamma_{i,j}^t e_t, & 1 \leq i, j \leq n, \\ [e_i, x_i] = e_i, & 1 \leq i \leq k, \\ [x_i, e_i] = (b_i - 1)e_i, & b_i \in \{0, 1\},\ 1 \leq i \leq k, \\ [e_i, x_j] = \alpha_{i,j} e_i, & k+1 \leq i \leq n,\ 1 \leq j \leq k, \\ [x_j, e_{i_t}] = \sum_{m=1}^{n_i} \beta_{j,t}^m e_{i_m}, & 1 \leq t \leq n_i,\ 1 \leq j \leq k, \end{cases} \quad (3.3)$$

*where omitted products are equal zero and $\alpha_{i,j}$ is the number of entries of a generator basis element $e_j$ involved in forming of non generator basis element $e_i$.*

This theorem implies the following corollary.

**Corollary 3.6.** *Let $R = N \oplus Q$ be a solvable Leibniz algebra such that $dimQ = dimN/N^2 = k$. Then $R$ admits a basis $\{e_1, e_2, \ldots, e_n, x_1, \ldots, x_k\}$ such that the vector space $N_i$ is invariant under the vector space $Q$ for all $i$.*

*Proof.* We will prove that $[N_i, x] \subseteq N_i$ and $[x, N_i] \subseteq N_i$, $1 \leq i \leq s-1$ for any $x \in Q$. By Theorem 3.5, there is a basis $\{e_1, e_2, \ldots, e_n, x_1, \ldots, x_k\}$ of the algebra $R$, in which the multiplication table has the form (3.3). Let $f \in N_i$, then $f$ is expressed by a linear combination in terms of basic elements $\{e_{i_1}, \ldots, e_{i_{n_i}}\} \subseteq N_i$, i.e. $f = \sum_{t=1}^{n_i} c_t e_{i_t}$. For any elements $x \in Q$, we have $x = \sum_{j=1}^{k} \mu_j x_j$. Then for $1 \leq t \leq n_i$ considering the following products:

$$[e_{i_t}, x] = [e_{i_t}, \sum_{j=1}^{k} \mu_j x_j] = \sum_{j=1}^{k} \mu_j [e_{i_t}, x_j] = (\sum_{j=1}^{k} \mu_j \alpha_{i_t, j}) e_{i_t}.$$

$$[x, e_{i_t}] = [\sum_{j=1}^{k} \mu_j x_j, e_{i_t}] = \sum_{j=1}^{k} \mu_j [x_j, e_{i_t}] = \sum_{m=1}^{n_i} \sum_{j=1}^{k} \mu_j \beta_{j,t}^m e_{i_m}.$$

Finally, by considering the following products, we conclude the proof of the corollary:

$$[f, x] = [\sum_{t=1}^{n_i} c_t e_{i_t}, \sum_{j=1}^{k} \mu_j x_j] = \sum_{t=1}^{n_i} \sum_{j=1}^{k} c_t \mu_j [e_{i_t}, x_j] = \sum_{t=1}^{n_i} \sum_{j=1}^{k} c_t \mu_j \alpha_{i_t, j} e_{i_t} \in N_i.$$

$$[x, f] = [\sum_{j=1}^{k} \mu_j x_j, \sum_{t=1}^{n_i} c_t e_{i_t}] = \sum_{t=1}^{n_i} \sum_{j=1}^{k} \mu_j c_t [x_j, e_{i_t}] = \sum_{m=1}^{n_i} \sum_{t=1}^{n_i} \sum_{j=1}^{k} c_t \mu_j \beta_{j,t}^m e_{i_m} \in N_i.$$

□

### 3.3. Solvable Leibniz algebras with a nilradical $\mathcal{L}(\alpha, \beta, \gamma)$ and the maximal codimensional is equal to one.

**Theorem 3.7.** *There is no solvable Leibniz algebra with the nilradical $\mathcal{L}(\alpha, \beta, \gamma)$ and the maximal dimension of the complementary space to the nilradical is equal to one.*

*Proof.* According to the condition, the maximal dimension of the complementary space of the solvable Leibniz algebra $R$ with a nilradical $\mathcal{L}(\alpha, \beta, \gamma)$ is equal to one. Using the Table 1, we get $a_{n-1} = 0$, $b_{n-1} = a_1$, $b_i = \alpha a_i$, $2 \leq i \leq n-3$ and $(\beta, \gamma) \neq (0, 0)$. Since $e_1, e_{n-1} \notin Ann_r(R)$, $e_2, e_3, \ldots, e_{n-2}, e_n \in Ann_r(R)$ and from Proposition 3.1 we have the following products in the algebra $R$:

$$[x, e_n] = 0, \ [x, e_1] = -e_1 + \sum_{t=2}^{n-2} c_{1,t} e_t + c_{1,n} e_n, \ [x, e_{n-1}] = \sum_{t=2}^{n-2} c_{n-1,t} e_t - e_{n-1} + c_{n-1,n} e_n.$$

Consider the following equality:

$$0 = [x, e_n] = [x, [e_{n-1}, e_1] - \alpha e_2] = [[x, e_{n-1}], e_1] - [[x, e_1], e_{n-1}] =$$

$$= [\sum_{t=2}^{n-2} c_{n-1,t} e_t - e_{n-1} + c_{n-1,n} e_n, e_1] - [-e_1 + \sum_{t=2}^{n-2} c_{1,t} e_t + c_{1,n} e_n, e_{n-1}] =$$

$$= \sum_{t=3}^{n-2} c_{n-1,t-1} e_t - (e_n + \alpha e_2) + \beta e_n - \sum_{t=3}^{n-2} c_{1,t-1} e_t.$$

From this we obtain:

$$\alpha = 0, \quad \beta = 1.$$

From the Table 1 it follows that $(\alpha, \beta) \neq (0, 1)$, i.e. this yields a contradiction. □

3.4. **Solvable Leibniz algebras with a nilradical $\mathcal{L}(\alpha, \beta, \gamma)$ and the maximal codimensional is equal to two.**

**Theorem 3.8.** *Let $R$ be a solvable Leibniz algebra with the nilradical $\mathcal{L}(\alpha, \beta, \gamma)$ and the maximal dimension of the complementary space to the nilradical be equal to two. Then $R$ is isomorphic to one of the following pairwise non-isomorphic algebras:*

$$R^1_{n+2}(0,\beta,0) : \begin{cases} [e_i, x] = ie_i, & 1 \leq i \leq n-2, \quad [e_n, x] = e_n, \quad [x, e_1] = -e_1, \quad [x, e_n] = \beta e_n, \\ [e_{n-1}, y] = e_{n-1}, \quad [e_n, y] = e_n, \quad [y, e_{n-1}] = \beta e_{n-1}, \quad [y, e_n] = \beta e_n, \quad \beta \in \{-1, 0\}, \end{cases}$$

$$R^2_{n+2}(0,1,1) : \begin{cases} [e_1, x] = e_1 - e_{n-1}, & [e_2, x] = 2e_2 - 2e_n, \quad [e_i, x] = ie_i, \quad 3 \leq i \leq n-2, \\ [x, e_1] = -e_1 + e_{n-1}, & [e_1, y] = e_{n-1}, \quad [e_2, y] = 2e_n, \quad [e_{n-1}, y] = e_{n-1}, \\ [e_n, y] = 2e_n, & [y, e_1] = -e_{n-1}, \quad [y, e_{n-1}] = -e_{n-1}, \end{cases}$$

$$R^3_{n+2}(1,-1,0) : \begin{cases} [e_1, x] = e_1 - e_{n-1}, & [e_i, x] = (i-1)e_i, \quad [e_n, x] = e_n, \quad [x, e_1] = -e_1 + e_{n-1}, \\ [x, e_n] = -e_n, & [e_1, y] = e_{n-1}, \quad [e_i, y] = e_i, \quad 2 \leq i \leq n-2, \\ [e_{n-1}, y] = e_{n-1}, & [e_n, y] = e_n, \quad [y, e_1] = -e_{n-1}, \quad [y, e_{n-1}] = -e_{n-1}, \\ [y, e_n] = -e_2 - e_n, \end{cases}$$

$$R^4_{n+2}(1,0,0) : \begin{cases} [e_1, x] = e_1 - e_{n-1}, & [e_2, x] = e_2 - e_n, \quad [e_i, x] = (i-1)e_i, \quad [e_n, x] = 2e_n, \\ [x, e_1] = -e_1 + e_{n-1}, & [e_1, y] = e_{n-1}, \quad [e_2, y] = e_2 + e_n, \quad [e_i, y] = e_i, \\ [e_{n-1}, y] = e_{n-1}, & & & 3 \leq i \leq n-2, \end{cases}$$

*where it is taken into account that each solvable algebra has its own multiplications of the nilradical and other products are zero.*

*Proof.* It is easy to see that $e_1$ and $e_{n-1}$ are the generator basis elements of the algebra $\mathcal{L}(\alpha, \beta, \gamma)$. So we have $\dim Q \leq 2$. By the hypothesis of the theorem, we need to investigate solvable Leibniz algebras with the dimension of the complementary subspace to the nilradical equal to the number of generators of the nilradical, i.e. $\dim Q = 2$. Let $\{x, y\}$ be a basis of the subspace $Q$. Then according to Theorem 3.5 and Corollary 3.6 we have the following brackets, i.e., the vector space $N_1$ is invariant under the vector space $Q$:

$$\begin{cases} [e_1, x] = e_1 + A_1 e_{n-1}, & [e_1, y] = B_1 e_{n-1}, \\ [e_{n-1}, x] = A_2 e_1, & [e_{n-1}, y] = B_2 e_1 + e_{n-1}, \\ [x, e_1] = \mu_{1,1} e_1 + \mu_{1,n-1} e_{n-1}, & [y, e_1] = \mu_{3,1} e_1 + \mu_{3,n-1} e_{n-1}, \\ [x, e_{n-1}] = \mu_{2,1} e_1 + \mu_{2,n-1} e_{n-1}, & [y, e_{n-1}] = \mu_{4,1} e_1 + \mu_{4,n-1} e_{n-1}, \\ [x, x] = [x, y] = [y, x] = [y, y] = 0. \end{cases}$$

From the Leibniz identity $\mathcal{L}I(e_1, x, y) = 0$, we get $B_1 = -A_1$. Taking into account that $\mathcal{R}_x$ and $\mathcal{R}_y$ are derivations of the algebra $\mathcal{G}(\alpha, \beta, \gamma)$ furthermore $e_1 \notin Ann_r(R)$ and $e_2, e_3, \ldots, e_{n-2} \in Ann_r(R)$, then multiplications in the solvable algebra $R$ have the following form:

$$\begin{cases} [e_1, x] = e_1 + A_1 e_{n-1}, & [e_1, y] = -A_1 e_{n-1}, \\ [e_2, x] = (2 + A_1\alpha)e_2 + A_1(1+\beta)e_n, & [e_2, y] = -A_1\alpha e_2 - A_1(1+\beta)e_n, \\ [e_i, x] = (i + A_1\alpha)e_i, & [e_i, y] = -A_1\alpha e_i, \quad 3 \leq i \leq n-2, \\ & [e_{n-1}, y] = e_{n-1}, \\ [e_n, x] = (1 + A_1\gamma - A_1\alpha(1+\beta))e_n, & [e_n, y] = (1 - A_1\gamma + A_1\alpha(1+\beta))e_n, \\ [x, e_1] = -e_1 + \mu_{1,n-1} e_{n-1}, & [y, e_1] = \mu_{3,n-1} e_{n-1}, \\ [x, e_{n-1}] = \mu_{2,n-1} e_{n-1}, & [y, e_{n-1}] = \mu_{4,n-1} e_{n-1}, \\ [x, e_n] = \sum_{t=2}^{n} \delta_{1,t} e_t, & [y, e_n] = \sum_{t=2}^{n} \delta_{2,t} e_t, \end{cases}$$

where $\alpha + A_1\alpha^2 = 0$, $\gamma(1 + A_1\gamma - A_1\alpha(1+\beta)) = 0$, $\gamma A_1 = \beta A_1(\gamma - \alpha(1+\beta))$.

Considering the Leibniz identity, we obtain the following restrictions on structure constants:

$$\begin{aligned}
\mathcal{L}I(x,x,e_{n-1}) &= 0, & &\Rightarrow & \mu_{2,n-1} &= 0, \\
\mathcal{L}I(x,e_1,y]) &= 0, & &\Rightarrow & \mu_{1,n-1} &= -A_1, \\
\mathcal{L}I(y,y,e_{n-1}) &= 0, & &\Rightarrow & \mu_{4,n-1}(1+\mu_{4,n-1}) &= 0, \\
\mathcal{L}I(y,e_1,y) &= 0, & &\Rightarrow & \mu_{3,n-1} &= -A_1\mu_{4,n-1}, \\
\mathcal{L}I(x,e_{n-1},e_1) &= 0, & &\Rightarrow & [x,e_n] &= (\beta + A_1\gamma)e_n, \\
\mathcal{L}I(y,e_{n-1},e_1) &= 0, & &\Rightarrow & [y,e_n] &= \mu_{4,n-1}\alpha e_2 + \mu_{4,n-1}(1+A_1\gamma)e_n, \\
\mathcal{L}I(x,y,e_n) &= 0, & &\Rightarrow & (\beta + A_1\gamma)(1 - A_1\gamma + A_1\alpha(1+\beta) + \mu_{4,n-1}(1+A_1\gamma)) &= 0, \\
\mathcal{L}I(x,e_{n-1},e_{n-1}) &= 0, & &\Rightarrow & \gamma(\beta + A_1\gamma) &= 0, \\
\mathcal{L}I(y,e_{n-1},e_{n-1}) &= 0, & &\Rightarrow & \gamma\mu_{4,n-1}\alpha &= 0, \quad \gamma\mu_{4,n-1}(1+A_1\gamma) = 0, \\
\mathcal{L}I(x,e_1,e_{n-1}) &= 0, & &\Rightarrow & \beta(1+\beta+A_1\gamma) + A_1\gamma &= 0, \\
\mathcal{L}I(y,e_1,e_{n-1}) &= 0, & &\Rightarrow & \mu_{4,n-1}\alpha(1+\beta) &= 0, \quad \mu_{4,n-1}(\beta(1+A_1\gamma)+A_1\gamma) = 0,
\end{aligned}$$

Thus, the table of multiplications of the algebra $R$ has form:

$$\begin{cases}
[e_1, x] = e_1 + A_1 e_{n-1}, & [e_1, y] = -A_1 e_{n-1}, \\
[e_2, x] = (2 + A_1\alpha)e_2 + A_1(1+\beta)e_n, & [e_2, y] = -A_1\alpha e_2 - A_1(1+\beta)e_n, \\
[e_i, x] = (i + A_1\alpha)e_i, & [e_i, y] = -A_1\alpha e_i, & 3 \le i \le n-2, \\
& [e_{n-1}, y] = e_{n-1}, \\
[e_n, x] = (1 + A_1\gamma - A_1\alpha(1+\beta))e_n, & [e_n, y] = (1 - A_1\gamma + A_1\alpha(1+\beta))e_n, \\
[x, e_1] = -e_1 - A_1 e_{n-1}, & [y, e_1] = -A_1\mu_{4,n-1}e_{n-1}, \\
& [y, e_{n-1}] = \mu_{4,n-1}e_{n-1}, \\
[x, e_n] = (\beta + A_1\gamma)e_n, & [y, e_n] = \mu_{4,n-1}\alpha e_2 + \mu_{4,n-1}(1+A_1\gamma)e_n,
\end{cases}$$

with

$$\begin{cases}
\alpha + A_1\alpha^2 = 0, & \gamma(1 + A_1\gamma - A_1\alpha(1+\beta)) = 0, \ \gamma A_1 = \beta A_1(\gamma - \alpha(1+\beta)), \\
\mu_{4,n-1}(1+\mu_{4,n-1}) = 0, & (\beta + A_1\gamma)(1 - A_1\gamma + A_1\alpha(1+\beta) + \mu_{4,n-1}(1+A_1\gamma)) = 0, \\
\gamma(\beta + A_1\gamma) = 0, & \gamma\mu_{4,n-1}\alpha = 0, \ \gamma\mu_{4,n-1}(1+A_1\gamma) = 0, \\
\mu_{4,n-1}\alpha(1+\beta) = 0, & \beta(1+\beta+A_1\gamma) + A_1\gamma = 0, \ \mu_{4,n-1}(1+\beta(1+A_1\gamma)+A_1\gamma) = 0.
\end{cases} \quad (3.4)$$

Using the Table 1, we have the following possible cases for $(\alpha, \beta, \gamma)$:

$$(\alpha, \beta, \gamma) \in \{(0, \beta, 0); (0, 1, 1); (1, -1, 0); (1, 0, 0)\}.$$

**Case 1.** Let $(\alpha, \beta, \gamma) = (0, \beta, 0)$. Then from (3.4) we conclude $\mu_{4,n-1} = \beta$ and by choosing $e'_1 = e_1 + A_1 e_{n-1}$, $e'_2 = e_2 + A_1(1+\beta)e_n$ we can assume $A_1 = 0$. Hence, we obtain the algebra $R^1_{n+2}(0, \beta, 0)$, where $\beta \in \{-1, 0\}$.

**Case 2.** Let $(\alpha, \beta, \gamma) = (0, 1, 1)$. Then $e_{n-1} \notin Ann_r(R)$ and $A_1 = \mu_{4,n-1} = -1$. Therefore, the algebra $R^2_{n+2}(0, 1, 1)$ is obtained.

**Case 3.** Let $(\alpha, \beta, \gamma) = (1, -1, 0)$. Then from restrictions (3.4) we get $A_1 = \mu_{4,n-1} = -1$. So, we obtain the algebra $R^3_{n+2}(1, -1, 0)$.

**Case 4.** Let $(\alpha, \beta, \gamma) = (1, 0, 0)$. Then using restrictions (3.4) we derive $A_1 = -1$ and $\mu_{4,n-1} = 0$. In this case we obtain the algebra $R^4_{n+2}(1, 0, 0)$. □

**Remark 3.9.** *The nilradical of the solvable algebra $R^1_{n+2}(0, \beta, 0)$ is $\mathcal{L}^{1,\beta}_n$. The classification of this solvable algebra is stated in paper [2]. Moreover, if we take $e'_1 = e_1 - e_{n-1}$, $e'_2 = e_2 - e_n$ in the algebra $R^2_{n+2}(0, 1, 1)$, then this algebra is isomorphic to the direct sum of solvable Leibniz algebras with null-filiform nilradical. Such Leibniz algebra was studied in the work [24, Corollary 3.4]:*

$$R^2_{n+2}(0,1,1) \cong (NF_{n-2} + \langle x \rangle) \oplus (NF_2 + \langle y \rangle):$$

$$\begin{cases} [e_i, e_1] = e_{i+1}, & 1 \le i \le n-3, \\ [e_i, x] = ie_i, & 1 \le i \le n-2, \\ [x, e_1] = -e_1, & \end{cases} \oplus \begin{cases} [e_{n-1}, e_{n-1}] = e_n, & [e_{n-1}, y] = e_{n-1}, \\ [e_n, y] = 2e_n, & [y, e_{n-1}] = -e_{n-1}. \end{cases}$$

3.5. **Solvable Leibniz algebras with a nilradical $\mathcal{G}(\alpha, \beta, \gamma)$ and the maximal codimensional is equal to one.**

**Theorem 3.10.** *Let $R$ be a solvable Leibniz algebra with the nilradical $\mathcal{G}(\alpha, \beta, \gamma)$ and the maximal dimension of the complementary space to the nilradical be equal to one. Then $R$ is isomorphic to one of the following pairwise non-isomorphic algebras:*

$$H_{n+1}^1(0,0,1) : \begin{cases} [e_1, x] = e_1, & [e_2, x] = 2e_2, & [e_i, x] = (i-2)e_i, \\ [x, e_1] = -e_1, & & [x, e_i] = -(i-2)e_i, & 3 \leq i \leq n. \end{cases}$$

$$H_{n+1}^2(1,2,0) : \begin{cases} [e_1, x] = e_1, & [e_2, x] = 2e_2, & [e_i, x] = (i-2)e_i, & [x, e_1] = -e_1, \\ [x, e_3] = -e_3, & [x, e_4] = -2e_4 + 2e_2, & [x, e_i] = -(i-2)e_i, & 3 \leq i \leq n. \end{cases}$$

$$H_{n+1}^3(1,0,\gamma) : \begin{cases} [e_1, x] = e_1, & [e_2, x] = 2e_2, & [e_i, x] = (i-2)e_i, \\ [x, e_1] = -e_1, & & [x, e_i] = -(i-2)e_i, & 3 \leq i \leq n. \end{cases}$$

$$H_{n+1}^4(1,-2,1) : \begin{cases} [e_1, x] = e_1, & [e_2, x] = 2e_2, & [e_i, x] = (i-2)e_i, & 3 \leq i \leq n, \\ [x, e_1] = -e_1, & [x, e_3] = -e_3, & [x, e_4] = -2e_4 - 2e_2, \\ & & [x, e_i] = -(i-2)e_i, & 5 \leq i \leq n. \end{cases}$$

$$H_{n+1}^5(1,4,2) : \begin{cases} [e_1, x] = e_1, & [e_2, x] = 2e_2, & [e_i, x] = (i-2)e_i, & 3 \leq i \leq n, \\ [x, e_1] = -e_1, & [x, e_3] = -e_3, & [x, e_4] = -2e_4 + 4e_2, \\ & & [x, e_i] = -(i-2)e_i, & 5 \leq i \leq n, \end{cases}$$

*where it is taken into account that each solvable algebra has its own multiplications of the nilradical and other products are zero.*

*Proof.* By condition theorem $R$ is a solvable Leibniz algebras with a codimension one nilradical $\mathcal{G}(\alpha, \beta, \gamma)$. Then using the above table, we get $b_3 = a_1$, $a_3 = 0$ and we have the following possible cases for $(\alpha, \beta, \gamma)$:

$$(\alpha, \beta, \gamma) \in \{(0, 0, 1); (1, 2, 0); (1, 0, \gamma); (1, -2, 1); (1, 4, 2)\}, \quad \gamma \neq 0.$$

From Proposition 3.2 we have the products in the algebra $R$:

$$\begin{cases} [e_1, x] = e_1 + a_2 e_2 + \sum_{t=4}^n a_t e_t, \\ [e_3, x] = b_2 e_2 + e_3 + \sum_{t=4}^n b_t e_t, \\ [e_2, x] = 2e_2, \\ [e_i, x] = (i-2)e_i + \sum_{t=i+1}^{n-1} b_{t-i+3} e_t + (b_{n-i+3} - (-1)^i a_{n-i+3}\alpha)e_n, & 4 \leq i \leq n-1, \\ [e_n, x] = (n-2)e_n, \\ [x, e_i] = \sum_{t=1}^n c_{n,i} e_t, & 1 \leq i \leq n, \\ [x, x] = \sum_{t=1}^n d_t e_t. \end{cases}$$

where $(1 - (-1)^n)a_{n-1}\alpha = 0$.

Applying the basis transformations in the following form:

$$e_1' = e_1 - a_2 e_2, \; e_2' = e_2, \; e_3' = e_3 - b_2 e_2 + \sum_{t=4}^n A_t e_t, \; e_i' = e_i + \sum_{t=i+1}^n A_{t-i+3} e_t, \; 4 \leq i \leq n,$$

with

$$A_4 = -b_4, \quad A_i = -\frac{1}{i-3}(b_i + \sum_{t=4}^{i-1} A_t b_{i-t+3}), \; 5 \leq i \leq n-1,$$

$$A_n = -\frac{1}{n-3}(b_n + \sum_{t=4}^n A_t(b_{n-t+3} - (-1)^t a_{n-t+3})),$$

we obtain $a_2 = b_2 = b_t = 0$ for $4 \leq i \leq n$.

Taking $e'_1 = e_1 + \beta a_4 e_2$, $x' = x + \sum_{t=4}^{n} a_t e_{t-1} - \frac{d_2}{2} e_2$, we can assume $d_2 = a_t = 0$ for $4 \leq t \leq n$.

It is easy to see that using products in the nilradical $\mathcal{G}(\alpha, \beta, \gamma)$, we have $e_1, e_3, \ldots, e_{n-1} \notin Ann_r(R)$ and $e_2 \in Ann_r(R)$. Thus, the table of multiplications of the algebra $R$ has the form:

$$\begin{cases} [e_1, x] = e_1, & [e_2, x] = 2e_2, & [e_i, x] = (i-2)e_i, & 3 \leq i \leq n, \\ [x, e_1] = -e_1 + c_{1,2} e_2 + c_{1,n} e_n, \\ [x, e_i] = -(i-2)e_i + c_{i,2} e_2 + c_{i,n} e_n, & 3 \leq i \leq n-1, & [x, e_n] = c_{n,2} e_2 + c_{n,n} e_n, & [x, x] = d_n e_n. \end{cases}$$

From the equalities $\mathcal{L}I(x, e_i, e_1) = \mathcal{L}I(x, x, x) = \mathcal{L}I(x, e_1, x) = \mathcal{L}I(x, e_3, x) = 0$ with $3 \leq i \leq n-1$, we derive the restrictions:

$$c_{4,2} = \beta, \ c_{4,n} = c_{i,2} = c_{i,n} = 0, \quad 5 \leq i \leq n-1, \quad c_{n,2} = 0, \ c_{n,n} = -(n-2),$$

$$d_n = c_{1,2} = c_{1,n} = c_{3,2} = c_{3,n} = 0.$$

Thus, the table of multiplications of the algebra $R$ has the form:

$$\begin{cases} [e_1, x] = e_1, & [e_2, x] = 2e_2, & [e_i, x] = (i-2)e_i, & 3 \leq i \leq n, \\ [x, e_1] = -e_1, & [x, e_3] = -e_3, & [x, e_4] = -2e_4 + \beta e_2, & [x, e_i] = -(i-2)e_i, & 5 \leq i \leq n. \end{cases}$$

Finally, we obtain solvable algebras $H^1_{n+1}(0, 0, 1)$, $H^2_{n+1}(1, 2, 0)$, $H^3_{n+1}(1, 0, \gamma)$, $H^4_{n+1}(1, -2, 1)$, $H^5_{n+1}(1, 4, 2)$ corresponding to the values of the parameter triples $(\alpha, \beta, \gamma)$, namely $(0, 0, 1)$, $(1, 2, 0)$, $(1, 0, \gamma)$, $(1, -2, 1)$, $(1, 4, 2)$, $\gamma \neq 0$.

□

**Remark 3.11.** *The nilradical of the solvable algebra $H^1_{n+1}(0, 0, 1)$ is $\mathcal{L}^3_n$. The classification of this solvable algebra is stated in paper [31].*

**3.6. Solvable Leibniz algebras with a nilradical $\mathcal{G}(\alpha, \beta, \gamma)$ and the maximal codimensional is equal to two.**

Let us give a classification of solvable Leibniz algebras with nilradical $\mathcal{G}(\alpha, \beta, \gamma)$ and two-dimensional complementary vector subspace to the nilradical.

**Theorem 3.12.** *Let $R$ be a solvable Leibniz algebra with the nilradical $\mathcal{G}(\alpha, \beta, \gamma)$ and the maximal dimension of the complementary space to the nilradical be equal to two. Then $R$ is isomorphic to one of the following pairwise non-isomorphic algebras:*

$$H^1_{n+2}(0, 0, 0) : \begin{cases} [e_1, x] = e_1, & [e_2, x] = 2e_2, & [e_i, x] = (i-3)e_i, & [x, e_1] = -e_1, \\ [x, e_i] = -(i-3)e_i, & [e_i, y] = e_i, & [y, e_i] = -e_i, & 3 \leq i \leq n. \end{cases}$$

$$H^2_{n+2}(0, 1, 0) : \begin{cases} [e_1, x] = e_1 - e_3, & [e_2, x] = 2e_2, & [e_i, x] = (i-3)e_i, & 4 \leq i \leq n, \\ [x, e_1] = -e_1 + e_3, & [x, e_4] = -e_4 + e_2, & [x, e_i] = -(i-3)e_i, & 5 \leq i \leq n, \\ [e_1, y] = e_3, & [e_i, y] = e_i, & 2 \leq i \leq n, \\ [y, e_1] = e_3, & [y, e_i] = -e_i, & 3 \leq i \leq n. \end{cases}$$

$$H^3_{n+2}(0, 2, 1) : \begin{cases} [e_1, x] = e_1 - e_3, & [e_4, x] = e_4 - e_2, & [e_i, x] = (i-3)e_i, & [x, e_1] = -e_1 + e_3, \\ [x, e_4] = -e_4 + e_2, & [x, e_i] = -(i-3)e_i, & [e_1, y] = e_3, & [e_2, y] = 2e_2, \\ [e_3, y] = e_3, & [e_4, y] = e_4 + e_2, & [e_i, y] = e_i, & [y, e_1] = -e_3, \\ [y, e_3] = -e_3, & [y, e_4] = -e_4 + e_2, & [y, e_i] = -e_i, & 5 \leq i \leq n. \end{cases}$$

$$H^4_{n+2}(1, 0, 0) : \begin{cases} [e_1, x] = e_1 - e_3, & [e_2, x] = 2e_2, & [e_i, x] = (i-3)e_i, & 4 \leq i \leq n, \\ [x, e_1] = -e_1 + e_3, & [x, e_i] = -(i-3)e_i, & & 4 \leq i \leq n, \\ [e_1, y] = -e_3, & [e_i, y] = e_i, & [e_n, y] = 2e_n, & 3 \leq i \leq n-1, \\ [y, e_1] = e_3, & [y, e_i] = -e_i, & [y, e_n] = -2e_n, & 3 \leq i \leq n-1. \end{cases}$$

$$H_{n+2}^5(1,1,0) : \begin{cases} [e_1, x] = e_1 - e_3, & [e_2, x] = e_2, & [e_i, x] = (i-3)e_i, & 4 \leq i \leq n, \\ [x, e_1] = -e_1 + e_3, & [x, e_4] = -e_4 + e_2, & [x, e_i] = -(i-3)e_i, & 5 \leq i \leq n, \\ [e_1, y] = e_3, & [e_i, y] = e_i, & [e_n, y] = 2e_n, & 2 \leq i \leq n-1, \\ [y, e_1] = -e_3, & [y, e_i] = -e_i, & [y, e_n] = -2e_n, & 3 \leq i \leq n-1, \end{cases}$$

$$H_{n+2}^6(1,2,1) : \begin{cases} [e_1, x] = e_1 - e_3, & [e_4, x] = e_4 - e_2, & [e_i, x] = (i-3)e_i, & 5 \leq i \leq n, \\ [x, e_1] = -e_1 + e_3, & [x, e_4] = -e_4 + e_2, & [x, e_i] = -(i-3)e_i, & 5 \leq i \leq n, \\ [e_1, y] = e_3, & [e_2, y] = 2e_2, & [e_3, y] = e_3, & [e_4, y] = e_4 + e_2, \\ [e_i, y] = e_i, & [e_n, y] = 2e_n, & & 5 \leq i \leq n-1, \\ [y, e_1] = -e_3, & [y, e_3] = -e_3, & [y, e_4] = -e_4 + e_2, & \\ [y, e_i] = -e_i, & [y, e_n] = -2e_n, & & 5 \leq i \leq n-1, \end{cases}$$

where it is taken into account that each solvable algebra has its own multiplications of the nilradical and other products are zero.

*Proof.* It is easy to see that $e_1$ and $e_3$ are the generator basis elements of the algebra $\mathcal{G}(\alpha, \beta, \gamma)$. So we have $\dim Q \leq 2$. By the hypothesis of the theorem, we need to investigate solvable Leibniz algebras with the dimension of the complementary subspace to the nilradical equal to the number of generators of the nilradical, i.e. $\dim Q = 2$. Let $\{x, y\}$ be a basis of the subspace $Q$. Then according to Theorem 3.5 and Corollary 3.6 we have the following brackets, i.e., the vector space $N_1$ is invariant under the vector space $Q$:

$$\begin{cases} [e_1, x] = e_1 + A_1 e_3, & [e_1, y] = B_1 e_3, \\ [e_3, x] = A_2 e_1, & [e_3, y] = B_2 e_1 + e_3, \\ [x, e_1] = \mu_{1,1} e_1 + \mu_{1,3} e_3, & [y, e_1] = \mu_{3,1} e_1 + \mu_{3,3} e_3, \\ [x, e_3] = \mu_{2,1} e_1 + \mu_{2,3} e_3, & [y, e_3] = \mu_{4,1} e_1 + \mu_{4,3} e_3, \\ [x, x] = [x, y] = [y, x] = [y, y] = 0. \end{cases}$$

Using the equality $\mathcal{L}I(e_1, x, y) = 0$, we deduce $B_1 = -A_1$. Taking into account that $\mathcal{R}_x$ and $\mathcal{R}_y$ are derivations of the algebra $\mathcal{G}(\alpha, \beta, \gamma)$ further $e_1, e_3, \ldots, e_{n-1} \notin \text{Ann}_r(R)$ and $e_2 \in \text{Ann}_r(R)$, then multiplications in the solvable algebra $R$ have the following form:

$$\begin{cases} [e_1, x] = e_1 + A_1 e_3, & [e_1, y] = -A_1 e_3, \\ [e_2, x] = (2 + A_1\beta)e_2, & [e_2, y] = -A_1\beta e_2, \\ & [e_3, y] = e_3, \\ [e_4, x] = e_4 + A_1\gamma e_2, & [e_4, y] = e_4 - A_1\gamma e_2, \\ [e_i, x] = (i-3)e_i, & [e_i, y] = e_i, & 5 \leq i \leq n-1, \\ [e_n, x] = (n-3-(-1)^n A_1\alpha)e_n, & [e_n, y] = (1 + (-1)^n A_1\alpha)e_n, \\ [x, e_1] = -e_1 - A_1 e_3, & [y, e_1] = A_1 e_3, \\ & [y, e_3] = -e_3, \\ [x, e_i] = -(i-3)e_i + \beta_{1,i} e_2 + \alpha_{1,i} e_n, & [y, e_i] = -e_i + \beta_{2,i} e_2 + \alpha_{2,i} e_n, & 4 \leq i \leq n-1, \\ [x, e_n] = \beta_{1,n} e_2 \alpha_{1,n} e_n, & [y, e_n] = \beta_{2,n} e_2 + \alpha_{2,n} e_n, \end{cases}$$

where non-written products are zero and

$$2A_1\gamma = \beta + A_1\beta^2, \quad \gamma(2 + A_1\beta) = 0, \quad \alpha(1 - (-1)^n A_1\alpha) = 0.$$

Using the Leibniz identity for the triples $\{x, e_i, e_1\}$ and $\{y, e_i, e_1\}$, $3 \leq i \leq n-1$, we conclude

$$\alpha_{1,i} = \alpha_{2,i} = 0, \quad 5 \leq i \leq n-1, \quad \beta_{1,i} = \beta_{2,i} = 0, \quad 5 \leq j \leq n,$$

$$\alpha_{1,n} = -(n-3-(-1)^n A_1\alpha), \quad \alpha_{2,n} = -(1+(-1)^n A_1\alpha),$$

$$[x, e_4] = -e_4 + (\beta + A_1\gamma)e_2, \quad [y, e_4] = -e_4 - A_1\gamma e_2.$$

Thus, the table of multiplications of the algebra $R$ has form:

$$\begin{cases} [e_1, x] = e_1 + A_1 e_3, & [e_1, y] = -A_1 e_3, \\ [e_2, x] = (2 + A_1\beta)e_2, & [e_2, y] = -A_1\beta e_2, \\ & [e_3, y] = e_3, \\ [e_4, x] = e_4 + A_1\gamma e_2, & [e_4, y] = e_4 - A_1\gamma e_2, \\ [e_i, x] = (i-3)e_i, & [e_i, y] = e_i, & 5 \leq i \leq n-1, \\ [e_n, x] = (n-3-(-1)^n A_1\alpha)e_n, & [e_n, y] = (1+(-1)^n A_1\alpha)e_n, \\ [x, e_1] = -e_1 - A_1 e_3, & [y, e_1] = A_1 e_3, \\ & [y, e_3] = -e_3, \\ [x, e_4] = -e_4 + (\beta + A_1\gamma)e_2, & [y, e_4] = -e_4 - A_1\gamma e_2, \\ [x, e_i] = -(i-3)e_i, & [y, e_i] = -e_i, & 5 \leq i \leq n-1, \\ [x, e_n] = -(n-3-(-1)^n A_1\alpha)e_n, & [y, e_n] = -(1+(-1)^n A_1\alpha)e_n. \end{cases}$$

where non-written products are zero and

$$2A_1\gamma = \beta + A_1\beta^2, \quad \gamma(2 + A_1\beta) = 0, \quad \alpha(1-(-1)^n A_1\alpha) = 0. \tag{3.5}$$

Using the Table 1, we have the following possible cases for $(\alpha, \beta, \gamma)$:

$$(\alpha, \beta, \gamma) \in \{(0,0,0); (0,1,0); (0,2,1); (1,0,0); (1,1,0); (1,2,1)\}.$$

**Case 1.** Let $(\alpha, \beta, \gamma) = (0,0,0)$. Then by choosing $e'_1 = e_1 + A_1 e_3$ we can assume $A_1 = 0$. Hence, we obtain the algebra $H^1_{n+2}(0,0,0)$.

**Case 2.** Let $(\alpha, \beta, \gamma) = (0,1,0)$. From the restrictions (3.5) we get $A_1 = -1$. So, we obtain the algebra $H^2_{n+2}(0,1,0)$.

**Case 3.** Let $(\alpha, \beta, \gamma) = (0,2,1)$. Then from (3.5) we conclude $A_1 = -1$ and obtain $H^3_{n+2}(0,2,1)$.

**Case 4.** Let $(\alpha, \beta, \gamma) = (1,0,0)$. Then $n$ is odd and in (3.5) we derive $A_1 = -1$. Therefore, the algebra $H^4_{n+2}(1,0,0)$ is obtained.

**Case 5.** Let $(\alpha, \beta, \gamma) = (1,1,0)$. Then $n$ is odd and using restrictions (3.5) we get $A_1 = -1$. Hence, we have the algebra $H^5_{n+2}(1,1,0)$.

**Case 6.** Let $(\alpha, \beta, \gamma) = (1,2,1)$, i.e $n$ is odd. Then from the above restrictions we obtain $A_1 = -1$ and the algebra $H^6_{n+2}(1,2,1)$.

□

**Remark 3.13.** *The nilradicals of the solvable algebras $H^1_{n+2}(0,0,0), H^2_{n+2}(0,1,0), H^3_{n+2}(0,2,1)$ and $H^4_{n+2}(1,0,0)$ given in Theorem 3.12 are $\mathcal{L}^1_n, \mathcal{L}^2_n, \mathcal{L}^4_n$ and $\mathcal{L}^5_n$, respectively. The classification of these solvable algebras is stated in papers [31], [32], [33] and [34], respectively.*

**Conclusion 1.** *Thus from the above table 1 and the obtained results it can be seen that the classifications of the solvable Leibniz algebras with the nilradical $\mathcal{L}^{3,-1}_n$, $\mathcal{L}^{3,0}_n$, $\mathcal{L}^{6,1}_n$, $\mathcal{L}^{8,2,1}_n$ and the dimension of complememtary space equals one have been remaining an open problem. For other 16 algebras the problem was solved.*

### 3.7. Solvable Leibniz algebras whose nilradical is a quasi-filiform Lie algebra.

In this subsection we describe solvable Leibniz algebras with the nilradical, naturally graded quasi-filiform Lie algebra, and the maximal dimension of complemented space of its nilradical. The whole class of complex Lie algebras $L$ having a naturally graded nilradical with characteristic sequence $C(L) = (n-2, 1, 1)$ is classified [21]. Here we find six families, two of them are decomposable, i.e., split into a direct sum of ideals: $\mathfrak{L}_{n,r}, \mathfrak{Q}_{n,r}, \mathfrak{T}_{n,n-3}, \mathfrak{T}_{n,n-4}, \mathfrak{L}_{n-1} \oplus \mathbb{C}$ and $\mathfrak{Q}_{2n} \oplus \mathbb{C}$, as well as there exist some special cases that appear only in low dimensions: $\mathfrak{E}_{7,3}, \mathfrak{E}^1_{9,5}, \mathfrak{E}^2_{9,5}$ and $\mathfrak{E}^3_{9,5}$.

We classify the solvable Leibniz algebras $R$ that have an indecomposable radical of arbitrary dimension.

**Theorem 3.14.** *Let $R$ be a solvable Leibniz algebra with the nilradical, natural graded quasi-filiform non-split Lie algebra, and the dimension of the complementary space to the nilradical be maximal. Then $R$ is a solvable Lie algebra.*

*Proof.* According to the condition, the nilradical of the solvable Leibniz $R$ is isomorphic to one of the following algebras:

$$\mathfrak{L}_{n,r}, \mathfrak{Q}_{n,r}, \mathfrak{T}_{n,n-3}, \mathfrak{T}_{n,n-4}, \mathfrak{E}_{7,3}, \mathfrak{E}^1_{9,5}, \mathfrak{E}^2_{9,5}, \mathfrak{E}^3_{9,5},$$

SOLVABLE LEIBNIZ ALGEBRAS WITH QUASI-FILIFORM NILRADICAL 15and the dimension of the complementary space is maximal.

Since the proof of the procedure repeats the same arguments that were presented earlier for each case, a detailed proof will be given only for the algebra $\mathfrak{L}_{n,r}$, the rest of the cases are completely analogous. Thus, we have

$$\mathfrak{L}_{n,r}(n \geq 5, \ r \text{ odd}, \ 3 \leq r \leq 2\left[\frac{n-1}{2}\right] - 1) : \begin{cases} [e_0, e_i] = e_{i+1}, & 1 \leq i \leq n-3, \\ [e_i, e_{r-i}] = (-1)^{i-1}e_{n-1}, & 1 \leq i \leq \frac{r-1}{2}, \end{cases}$$

where $\{e_0, e_1, \ldots, e_{n-1}\}$ is a basis of the algebra $\mathfrak{L}_{n,r}$.

The derivations of the algebra $\mathfrak{L}_{n,r}(n > 5)$ has the following form [20].

$$\begin{aligned}
&t_0(e_0) = e_0, \quad t_0(e_i) = (i-1)e_i, \quad 2 \leq i \leq n-2, \quad t_0(e_{n-1}) = (r-2)e_{n-1}; \\
&t_1(e_0) = e_1, \quad t_1(e_r) = e_{n-1}; \\
&t_2(e_i) = e_i, \quad 1 \leq i \leq n-2, \quad t_2(e_{n-1}) = 2e_{n-1}; \\
&h_k(e_i) = e_{k+i}, \quad 1 \leq i \leq n-2-k, \quad 3 \leq k \leq n-3, \quad k \text{ odd if } k \leq r-3; \\
&g_1(e_0) = e_r, \quad g_2(e_0) = e_{n-1}.
\end{aligned}$$

The shape of the derivations further shows that there are two independent non-nilpotent derivations $t_0$ and $t_2$. Let $\mathcal{R}_x$ and $\mathcal{R}_y$ denote the two nil-independent derivations of $\mathfrak{L}_{n,r}$ and let $\{x, y\}$ be a basis of the subspace $Q$. Then according to the derivations of $\mathfrak{L}_{n,r}$ and using $e_0, e_1, \ldots, e_{n-3} \notin Ann_r(R)$, we have the following brackets:

$$\begin{cases}
[e_0, x] = e_0 + a_1 e_1 + a_r e_r + a_{n-1} e_{n-1}, & [e_0, y] = b_1 e_1 + b_r e_r + b_{n-1} e_{n-1}, \\
[e_i, x] = (i-1)e_i + \sum_{t=i+1}^{n-1}(*)e_t, & [e_i, y] = e_i + \sum_{t=i+1}^{n-1}(*)e_t, \\
[e_{n-1}, x] = (r-2)e_{n-1}, & [e_{n-1}, y] = 2e_{n-1}, \\
[x, e_0] = -e_0 - a_1 e_1 - a_r e_r + c_{n-1}e_{n-1}, \\
[x, e_i] = -(i-1)e_i + \sum_{t=i+1}^{n-1}(*)e_t, & 1 \leq i \leq n-2.
\end{cases}$$

Using the equalities $\mathcal{L}I(x, e_0, e_{n-3}) = \mathcal{L}I(x, e_1, e_{r-1}) = 0$, we deduce $e_{n-2}, e_{n-1} \notin Ann_r(R)$. Since the ideal $I = \{[x,x] \mid x \in R\}$ is contained in $Ann_r(R)$, then we have $I = \{0\}$. Thus, we have shown that $R$ is a solvable Lie algebra. These solvable Lie algebras are studied in [5].

It should be noted that it is convenient to consider solvable Leibniz algebras with the nilradical $\mathfrak{L}_{5,3}$ separately from the general one. The reason is that the derivations of the algebra $\mathfrak{L}_{5,3}$ are different from those of the higher dimensional algebras. The classification of this solvable algebra is given in the article [30]. □

**Remark 3.15.** *The classification of the solvable Leibniz algebra with the nilradical $\mathfrak{L}_{n-1} \oplus \mathbb{C}$ or $\mathfrak{Q}_{2n} \oplus \mathbb{C}$ and the complementary space to the nilradical with a maximal dimension is stated in paper [22].*

In conclusion we formulate

**Conjecture.** Let $R$ be a solvable Leibniz algebra with the nilradical non-split Lie algebras and the dimension of the complementary space to the nilradical be maximal. Then $R$ is a solvable Lie algebra.